\documentclass[12pt,a4paper,leqno]{article}
\usepackage[english]{babel}
\usepackage{amscd,amsmath,amssymb,amsfonts}

\usepackage{epsf,pstricks,epsfig}

\usepackage{epic,eepic}
\usepackage{graphicx}

\def\leftmapsto{\leftarrow\hspace{-1,7mm}\raisebox{0,6pt}{$\shortmid$}}

\def\Spec{\mbox{Spec}}

\input xy
\xyoption{all}

\def\Spec{\mathrm{Spec}}
\def\dim{\mathrm{dim}}
\def\min{\mathrm{min}}
\def\max{\mathrm{max}}
\def\ch{\mathrm{char}}
\def\q{{\mathbb Q}}
\def\ql{{\mathbb Q_\ell}}
\def\qp{{\mathbb Q_p}}

\def\qzl{\mathbb Q_\ell/\mathbb Z_\ell}

\def\zl{{\mathbb Z_\ell}}

\def\z{\mathbb Z}

\def\qzl{\mathbb Q_\ell/\mathbb Z_\ell}
\def\ker{\mathrm{ker}}
\def\coker{\mathrm{coker}}

\def\Cb{\overline{C}}
\def\Hb{\overline{H}}
\def\ldot{\vspace{-1mm}\cdot}

\def\ch{\mathrm{char}}

\numberwithin{equation}{section}
   \newtheorem{thm}{Theorem}[section]
   
   \newtheorem{prop}[thm]{Proposition}
   \newtheorem{cor}[thm]{Corollary}
   \newtheorem{conj}[thm]{Conjecture}
   \newtheorem{ques}[thm]{Question}

\newcommand{\cqfd}

\begin{document}

\title
{Weights in Arithmetic Geometry}
\author{Uwe Jannsen}
\date{}
\maketitle
\thispagestyle{empty}

\begin{abstract}
The concept of weights on the cohomology of algebraic varieties
was initiated by fundamental ideas and work of A. Grothendieck and P.
Deligne. It is deeply connected with the concept of motives and appeared
first on the singular cohomology as the weights of (possibly mixed) Hodge
structures and on the etale cohomology as the weights of eigenvalues of
Frobenius. But weights also appear on algebraic fundamental groups and
in $p$-adic Hodge theory, where they become only visible after applying the
comparison functors of Fontaine. After rehearsing various versions of
weights, we explain some more recent applications of weights,
e.g., to Hasse principles and the computation of motivic cohomology,
and discuss some open questions.
\end{abstract}

\bigskip\bigskip
The theory of weights is already explained by Deligne in his talks at the
ICM at Nice \cite{De1} and the ICM at Vancouver \cite{De4},
in a magnificently concise and clear way, and every reader
is urged to read his account before starting with this paper. In addition
Deligne contributed the basic substance to this theory, by proving the Weil
conjectures in a very general way and establishing the theory of mixed
Hodge structures. So the modest aim of this article is just to give a certain update
of the results, and to discuss some applications of weights in arithmetic geometry.
For this we have of course to rehearse some of the theory of weights, at least
so far that our methods and results become clear.
A certain emphasis is on the fact that weights are intimately linked with
resolution of singularities, although this is not so clear from Deligne's
proof of the Weil conjecture. But resolution is, at least at present state, indispensable
in Hodge theory, and should also play a major role in establishing good
$p$-adic theories.

\bigskip
In the following, the word variety will mean a separated scheme of finite type
over a field $K$.

\section{Weights in Hodge theory}

\medskip
A \emph{pure} $\mathbb Q$-\emph{Hodge structure of weight} $n$ is a finite dimensional $\mathbb Q$-vector space $V$ together with a decomposition
$$
V_\mathbb C\; :=\; V\otimes_\mathbb Q\mathbb C\; =\; \mathop\oplus\limits_{p+q=n}H^{p,q}
$$
into $\mathbb C$-vector spaces such that $\overline{H^{p,q}}=H^{q,p}$ where $\overline{H^{p,q}}=(id\otimes\sigma)H^{p,q}$ for the complex conjugation $\sigma$ on $\mathbb C$. This is exactly the structure one gets on the $\mathbb Q$-cohomology $H^n(X,\mathbb Q)$ of a smooth projective complex variety $X$, the isomorphism with the de Rham cohomology
$$
H^n(X,\mathbb Q)\otimes_\mathbb Q\mathbb C\; \cong\; H^n_{dR}(X)
$$
and the famous Hodge decomposition
$$
H^n_{dR}(X)\; =\; \mathop\oplus_{p+q=n}H^{p,q}
$$
into the spaces of harmonic $(p,q)$-forms. A morphism $V \rightarrow V'$ of pure Hodge structures
is a $\q$-linear map respecting the bigrading after tensoring with $\mathbb C$. There is a natural tensor
product of Hodge structures, by defining $(V\otimes W)_{\mathbb C}^{p,q}$ as the sum over all
$V_{\mathbb C}^{i,j}\otimes W_{\mathbb C}^{k,l}$ with $i+k=p$ and $j+l=q$. One has the obvious

\medskip
\textbf{Fact 1} If $V$ and $W$ are pure Hodge structures of weights $m\neq n$, then every morphism
$\varphi: V \rightarrow W$ of Hodge structures is trivial.

\medskip
P. Deligne established the theory of \emph{mixed Hodge structures} for arbitrary complex varieties. These are finite dimensional $\mathbb Q$-vector spaces $V$ with an ascending filtration $W.V$ by $\mathbb Q$-vector spaces $(W_{n-1}V\subseteq W_n V)$ and a descending filtration $F^\cdot V_\mathbb C$ by $\mathbb C$-vector spaces such that

$$
Gr^W_n V\; :=\; W_n V/W_{n-1}V
$$
gets a pure Hodge structure of weight $n$ via $H^{p,q}_n=F^p\cap \overline{F^q}Gr^W_nV_{\mathbb C}$, for the induced filtration $F^\cdot $ on
$Gr^W_n V_{\mathbb C}$.
A morphism $V \rightarrow V'$ of mixed Hodge structures is a $\q$-linear map respecting both filtrations.
It is a non-trivial fact that this is an abelian category.
The pure Hodge structures form a full subcategory of the mixed Hodge structures.

\medskip
For example a smooth quasi-projective complex variety $U$ gets a mixed Hodge structure on its cohomology as follows. By Hironaka's resolution of singularities, there is a \emph{good compactification} of $U$, i.e., an open embedding $U\subset X$ into a smooth projective variety such that $Y=X-U$ is a simple normal crossings divisor, i.e., has smooth (projective) irreducible components $Y_1,\ldots , Y_N$ such that each $p$-fold intersection $Y_{i_1}\cap \ldots \cap Y_{i_p}$ (with $1\leq i_1<i_2<\ldots <i_p\leq N$) is smooth of pure codimension $p$. This gives rise to a combinatorial spectral sequence
\begin{equation}\label{eq.wt.ss.hodge}
E_2^{p,q}\; =\; H^p(Y^{[q]},\mathbb Q(-q))\quad\Rightarrow\quad H^{p+q}(U,\mathbb Q)
\end{equation}
where
$$
Y^{[q]}\; = \; \coprod\limits_{i_1<\ldots <i_q}Y_{i_1}\cap \ldots \cap Y_{i_q}
$$
is the disjoint union of all $q$-fold intersections of the $Y_i$. Here $\mathbb Q(-q)=\mathbb Q(2\pi i)^{-q}$, but moreover, one regards
$$
H^p(Y^{[q]},\mathbb Q(-q))\; \cong \; H^p(Y^{[q]},\mathbb Q)\otimes_\mathbb Q\mathbb Q(-q)
$$
as a pure Hodge structure of weight $p+2q$, via defining the \emph{Tate Hodge structure} $\mathbb Q(1)=\mathbb Q2\pi i$ as pure of weight $-2$ and Hodge type $(-1,-1)$, so that $\mathbb Q(m):=\mathbb Q(1)^{\otimes m}$, for $m\in\mathbb Z$, becomes pure of weight $-2m$ and type $(-m,-m)$ for $m\in\mathbb Z$, and the $m$-th Tate twist $V(m):=V\otimes \mathbb Q(m)$ of a pure Hodge structure $V$ of weight $n$ has weight $n-2m$. With this convention Deligne showed that all differentials $d_r^{p,q}$ in this spectral sequence are morphisms of pure Hodge structures, and that there is a mixed Hodge structure on $H^n(U,\mathbb Q)$ such that the weight filtration comes, up to certain shift, from the above spectral sequence: If $\widetilde W_\cdot$ is the ascending filtration on $H^n(U,\mathbb Q)$ with $\widetilde W_q/\widetilde W_{q-1}\cong E_\infty^{n-q,q}$ then the weight filtration $W_q$ is given by $W_q=\widetilde W_{q-n}$.

\medskip
A remarkable consequence, observed by Deligne, is that the spectral sequence degenerates at $E_3$, i.e., the differentials $d_r^{p,q}$ vanish for $r\geq 3$. In fact they go from a subquotient of $E^{p,q}_2=H^p(Y^{[q]},\mathbb Q(-q))$, of weight $p-2q$, to a subquotient of
$$
E_2^{p+r,q-r+1}\; =\; H^{p+r}(Y^{[q-r+1]}, \mathbb Q(-q+r-1))\,,
$$
which is of weight $p+r+2q-2r+2=p+2q-r+2\neq p+2q$ for $r\geq 3$, and we can use Fact 1. I do not know of any proof not using weights.

\medskip
A second remarkable consequence is that the $E_3=E_\infty$-terms do not depend on the choice of the good compactification. For this one has another proof using the so-called strong factorization in the theory of resolution of singularities.

\medskip
A third remarkable consequence is that the intrinsic structure of $H^n(U,\mathbb Q)$ in some sense ``sees'' the good compactification, at least the homology of the complexes of $E_2$-terms. For example, if $S$ is a surface with one isolated singular point $P$ which can be resolved by blowing up $P$, then the cohomology of $U=S-P$ sees the created exceptional divisor.

\medskip
It is not known if the weight filtration on $H^n(U,\mathbb Q)$ can be obtained by some other process, e.g., one that is intrinsic on $U$, or one which uses other types of compactification, e.g., some which appear in the minimal model program.

\medskip
Thus the theory of mixed Hodge structures and their weights depends on resolution of singularities.

\vspace{0,5cm}
\section{Weights in $\ell$-adic cohomology}

\medskip
 Let $K$ be an arbitrary field with separable closure $\overline K$, and let $\ell$ be a prime invertible in $K$.
Let $X/K$ be a smooth projective variety and let $\overline X=X\times_K\overline K$ be the base-change of $X$,
which is a smooth projective variety over $\overline{K}$.
Then the $\ell$-adic \'{e}tale cohomology
\begin{equation}\label{eq.l-adic.coh}
H^n(\overline X,\mathbb Q_\ell)
\end{equation}
is a finite-dimensional $\mathbb Q_\ell$-vector space with a continuous action of the absolute Galois group $G_K=Gal(\overline K/K)$, via functoriality of \'{e}tale cohomology: for $\sigma\in G_K$, $id\times\sigma$ acts on $X\times_K\overline K$ and induces $(id\times\sigma)^\ast$ on $H^n(\overline X,\mathbb Q_\ell)$.

\medskip
Now consider the case that $K = \mathbb F_q$ is a finite field with $q$ elements.
Then, as conjectured by A. Weil \cite{Wei} and reformulated and proved by P. Deligne \cite{De2},\cite{De3},
the Galois representation \eqref{eq.l-adic.coh} is \emph{pure of weight} $n$. This means: If $F\in G_K$ is a \emph{geometric} Frobenius automorphism -- the inverse of the \emph{arithmetic} one sending $x\in\overline K$ to $x^q$ -- then the eigenvalues $\alpha$ of $F$ (precisely: $F^\ast$) on $H^n(\overline X,\mathbb Q_\ell)$ are \emph{pure of weight} $n$, i.e., they are algebraic numbers (the characteristic polynomial of $F$ lies in $\mathbb Q[T]$ rather than in $\mathbb Q_\ell[T]$), and one has
\begin{equation}\label{eq.purity}
|\sigma(\alpha)|\; =\; q^{\frac n 2}
\end{equation}
for any field embedding $\sigma:\mathbb Q(\alpha)\hookrightarrow \mathbb C$. We mention in passing that
Deligne also proved independence of $\ell$: The characteristic polynomial does not depend on
$\ell$ ($\neq p = \ch(K)$).

\medskip
If $K$ is finitely generated (over its prime field), then $H^n(\overline X,\mathbb Q_\ell)$ is again pure of weight $n$, in the following sense: By assumption there is an integral scheme $S$ of finite type over $\mathbb Z$ with function field $K$. After possibly shrinking $S$ we may assume there is a smooth projective scheme $\mathcal X$ over $S$ such that $X=\mathcal X\times_SK$ (the fibre of $\pi:\mathcal X\rightarrow S$ over the generic
point $\eta = \Spec(K)$ of $S$). Then for every closed point $s\in S$ the theory of smooth and proper base change gives an isomorphism
\begin{equation}\label{eq.bc.l-adic}
b_{\eta,s}\; :\; H^n(X \times_K\overline K,\mathbb Q_\ell)\; \stackrel\sim\longrightarrow \; H^n(\mathcal X_s\times_s\overline{k(s)},\mathbb Q_\ell)
\end{equation}
where $\mathcal X_s=\pi^{-1}(s)$ is the fibre of $\mathcal X$ over $s$, $k(s)$ is the residue field of $s$, and $\overline{k(s)}$ is its separable closure. The above isomorphism is compatible with the Galois actions, in the following sense. There is a decomposition group $G_s\subseteq G_K$ with an epimorphism $\rho:G_s\twoheadrightarrow G_{k(s)}$ such that $b_{\eta,s}$ is compatible with the actions of $G_s$ and $G_{k(s)}$, respectively. So the inertia group $I_s=\ker(\rho)$ acts trivially on $H^n(\overline X,\mathbb Q_\ell)$ and via $\rho$, $b_{\eta,s}$ is a $G_{k(s)}$-morphism. Note that $k(s)$ is a finite field so that it makes sense to speak of a pure $\mathbb Q_\ell$-representation for $G_{k(s)}$.

\vspace{0,5cm}
A general $G_K$-representation $V$ is now called pure of weight $n$, if there is an $S$ above with function field $K$ such that for all $s\in S$ the inertia group $I_s\subseteq G_s$ acts trivially and the obtained $G_{k(s)}=G_s/I_s$-representation $V$ is pure of weight $n$.
Obviously we have:

\medskip
\textbf{Fact 2} If $V$ and $W$ are pure $\mathbb Q_\ell$-$G_K$-representations of weights $n\neq m$, then every $G_K$-homomorphism $\varphi:V\rightarrow W$ is zero.

\medskip
Again Deligne extended this to arbitrary varieties $X$ over $K$. First of all he introduced the concept of a mixed $\mathbb Q_\ell$-$G_k$-representation $V$; this simply means that there is a filtration $W'_n$ on $V$ such that each $Gr^{W'}_n V$ is pure (of some weight). For $K$ (finitely generated) of positive characteristic he showed that then one also has an ascending filtration $W_n$ on $V$ such that $Gr^W_nV$ is pure of weight $n$. A filtration with this property is unique (by Fact 2) and is called the weight filtration, but it could be different from $W'_n$ (look at the sum of two pure representations). Over a field $K$ of characteristic $0$ one can again use resolution of singularities to produce a weight filtration. E.g., for a smooth quasi-projective variety $U/K$ one can again choose a good compactification $U\subseteq X\supseteq Y=X-U$ to obtain a spectral sequence
\begin{equation}\label{eq.wt.ss.l-adic}
E_2^{p,q}\; =\; H^p(\overline{Y^{[q]}}, \mathbb Q_\ell(-q))\; \Rightarrow \; H^{p+q}(\overline U,\mathbb Q_\ell)
\end{equation}
completely analogous to \eqref{eq.wt.ss.hodge}. Here $\mathbb Q_\ell(m)=\mathbb Q_\ell(1)^{\otimes m}$, where $\mathbb Q_\ell(1)=\mathbb Q_\ell\otimes_{\mathbb Z_\ell}\mathbb Z_\ell (1)$ for $\mathbb Z_\ell(1)=\lim\limits_{\leftarrow n} \mu_{\ell^n}$, the inverse limit of the $G_K$-modules $\mu_{\ell^n}\subseteq \overline K$ of $\ell^n$-th roots of unity. Hence $\mathbb Z_\ell(1)$ is non-canonically isomorphic to $\mathbb Z_\ell$, and $G_K$ acts via the cyclotomic character. This definition makes sense over any field of characteristic $\neq \ell$. If $F$ is a finite field, the arithmetic Frobenius acts on $\mathbb Z_\ell(1)$ as multiplication by $q=|F|$. Therefore for any finitely generated $K$ with $\ell\neq \ch(K)$, $\mathbb Q_\ell(1)$ is pure of weight $-2$. Hence
$$
H^p(\overline{Y^{[q]}}, \mathbb Q_\ell (-q))\; \cong \; H(\overline{Y^{[q]}},\mathbb Q_\ell)(-q)
$$
is pure of weight $p+2q$. Here $V(m) := V\otimes \ql(m)$ is the $m$-th Tate twist of a $\ql$-$G_K$-representation. Like in \eqref{eq.wt.ss.hodge}, one obtains a weight filtration on $H^n(\overline U,\mathbb Q_\ell)$ as $W_i=\widetilde W_{i-n}$, where $\widetilde W_\cdot$ is the canonical ascending filtration associated to \eqref{eq.wt.ss.l-adic}.

\vspace{0,5cm}
Similar remarks as in the Hodge theoretic setting apply -- except that now it is clear that the weight filtration -- once it exists -- is intrinsic.

\vspace{0,5cm}
A noteworthy fact is that -- in contrast to the case of positive characteristic -- there exist non-trivial extensions.
$$
0\; \rightarrow\; V_1\; \rightarrow\; V\; \rightarrow \; V_2\;\rightarrow \; 0
$$
of $\ell$-adic representations over fields $K$ of characteristic $0$ such that $V_i$ is pure of weight $n_i$,
with $n_1>n_2$.
%We will come back to this phenomenon later.

\vspace{0,5cm}
\section{Weights in $p$-adic cohomology}

\medskip
Let $k$ be a perfect field of characteristic $p>0$, let $W=W(k)$ be the ring of Witt vectors over $k$, and let $K_0=$ Frac$(W)$ be
the fraction field of $W$. If $X$ is a smooth projective variety over $k$, then the crystalline cohomology
$$
H^n(X/K_0) := \; H^n (X/W)\otimes_W K_0
$$
is a finite-dimensional vector space over $K_0$.

\medskip
We shall need the following variant cohomology. If $W\Omega^\cdot_X$ is the de Rham-Witt complex of $X$,
see \cite{Il}, which is a pro-complex formed by the complexes $W_m\Omega^\cdot_X$ for all $m\in \mathbb N$, then one has a natural isomorphism
$$
H^n(X,W\Omega^\cdot _X)\; \stackrel\sim\longrightarrow  \; H^n(X/W)\, ,
$$
where the left hand side is \'{e}tale (hyper)cohomology. Let $W_m\Omega^r_{X,\log}\subseteq W_m\Omega^r_X$ be the logarithmic part of the de Rham-Witt sheaf -- it is  \'{e}tale locally generated by sections of the form $d \hat{x_1}/\hat{x_1}\wedge \ldots \wedge d\hat{x_r}/\hat{x_r}$, where $x_i\in \mathcal O_X^\times$ and $\hat{x_i}$ is a Teichm\"uller lift in $W_m\mathcal O_X$. Then there is an exact sequence of pro-sheaves (\cite{Il} I 5.7.2),
where $Fr$ is the Frobenius operator on the de Rham-Witt complex
$$
0\; \rightarrow\; W\Omega^r_{X,\log}\; \rightarrow\; W\Omega^r_X\; \stackrel{1-Fr}{\longrightarrow}\; W\Omega^r_X\; \rightarrow \; 0\, .
$$
If $k$ is algebraically closed, then it induces exact sequences (\cite{Mi1} 1.15)
\begin{equation}\label{eq.ex.seq.qp}
0\; \rightarrow\; H^n(X,\mathbb Q_p(r))\; \rightarrow H^n_{cris}(X/W)_{\mathbb Q_p}\; \stackrel{F-p^r}\longrightarrow\; H^n_{cris}(X/W)_{\mathbb Q_p}\; \rightarrow \; 0
\end{equation}
where $F$ is the (morphism induced by the) Frobenius endomorphism of $X$. Here we use the following notation by J. Milne (loc. cit.):
$$
\begin{array}{rcl}
\mathbb Z/p^m\mathbb Z(r) & = & W_m\Omega^r_{X,\log}[-r]\vspace{2mm}\\
H^n(X,\mathbb Z_p(r)) & = & \lim\limits_{\leftarrow\;m}
%\lim\limits_{\stackrel{\textstyle\leftarrow}{m}}
H^n(X,\mathbb Z/p^m\mathbb Z(r))\vspace{2mm}\\
H^n(X,\mathbb Q_p(r)) & = & H^n(X,\mathbb Z_p(r))\otimes_{\mathbb Z_p}\mathbb Q_p\, .
\end{array}
$$

\medskip
Now we discuss weights. First let $k$ be a finite field, and let $F:X\rightarrow X$ be the (geometric) Frobenius relative to $k$, i.e., the $k$-morphism which sends a local section $f$ of $\mathcal O_X$ to $f^q$, $q=|k|$.
Then by a simple but ingenious argument N. Katz and W. Messing \cite{KM} deduced from
Deligne's proof of the Weil conjectures that
$H^n(X/K_0)$ is pure of weight $n$ in the sense that the eigenvalues $\alpha$
of the endomorphism $F^\ast$ induced on it are algebraic numbers with absolute value
$|\sigma(\alpha)|=q^{\frac n 2}$ for every embedding $\sigma:\mathbb Q(\alpha)\hookrightarrow \mathbb C$.
Moreover the characteristic polynomial of $F^\ast$ is the same as the one obtained on
the $n$-th $\ql$-cohomology.

\medskip
By a lemma in (semi-) linear algebra, \cite{Mi1}, 5.1, the exact sequence \eqref{eq.ex.seq.qp} for $\overline X=X\times_k\overline k$ implies that the (linear!) action of the geometric Frobenius $F\in G_k$ on $H^n(\overline X,\mathbb Q_p(r))$ is pure of weight $n-2r$, in perfect analogy to the $\ell$-adic case. More precisely, the eigenvalues of $F$ on $H^n(\overline X,\mathbb Q_p(r))$ are those eigenvalues for $H^n(\overline X/K_0)$ which have slope $p^r$.

\medskip
Now let $K$ be a finitely generated field of characteristic $p>0$, and let $S$ (of finite type over $\mathbb F_p)$ and $\mathcal X\rightarrow S$ (smooth and proper) be as above. Then M. Gros and N. Suwa (\cite{GrSu} Th. 2.1) established a base change isomorphism
\begin{equation}\label{eq.bc.p-adic}
b_{\eta,s}: \,H^n(X\times_K\overline{K},\mathbb Q_p(r))\; \mathop{\longrightarrow}\limits^{\sim}\; H^n(\mathcal X_s\times_{k(s)}\overline{k(s)}, \mathbb Q_p(r))
\end{equation}
for every closed $s$ in a non-empty open subset $U \subseteq S$, which is equivariant for the decomposition group at $s$.
By this the left hand side is a pure $\qp$-$G_K$-representation in exactly the same sense as in the $\ell$-adic setting.

\vspace{0,5cm}
The case of a general variety is more difficult for $p$-adic cohomology. The crystalline cohomology does not behave well for singular or nonproper varieties; in particular it is not in general finite-dimensional. A good finite-dimensional theory is given by the rigid analytic cohomology
(this follows from de Jong's resolution of singularities \cite{deJ}, see P. Berthelot \cite{Be}), and it coincides with the crystalline cohomology for smooth proper varieties. But only quite recently purity and a theory of weights has been studied thoroughly in this context \cite{NS}, \cite{Nak},
again using \cite{deJ}.
Below I am rather interested in the cohomology $H^n(\overline X,\mathbb Q_p(r))$, but one does not have purity for general $r$, see M. Gros \cite{Gr}, and no theory of weights for general varieties either. Fortunately there is a good situation for $r=\dim(X)$ (loc. cit.)
which we will use for our applications.

\vspace{0,5cm}
\section{Weights for $\ell$-adic cohomology over local fields}

\medskip
Let $K$ be a non-archimedean local field, i.e., a complete discrete valuation field with finite residue field $k$. Let $X$ be a smooth projective variety over $K$ and let $\ell$ be a prime, $\ell\neq p=\ch(k)$. If $X$ has good reduction, i.e., a smooth proper model $\mathcal X\rightarrow\Spec(\mathcal O_K)$ over the ring of integers $\mathcal O_K$ (i.e., the discrete valuation ring) of $K$, then as before the base change isomorphism
\begin{equation}\label{eq.spec.l-adic}
H^n(\overline X,\mathbb Q_\ell)\; \mathop{\longrightarrow}\limits^{\sim}\; H^n(\overline{\mathcal X_k},\mathbb Q_\ell)
\end{equation}
shows that the inertia group $I\subset G_K$ acts trivially on $H^n(\overline X,\mathbb Q_\ell)$ and that $H^n(\overline X,\mathbb Q_\ell)$ corresponds to a pure $\mathbb Q_\ell$-representation of weight $n$ of $G_K/I\cong G_k$.

\vspace{0,5cm}
In general it is a theorem of Grothendieck that, after possibly passing to a finite separable extension of $K$, the ramification group $P\subset I$ acts trivially and the pro-cyclec group $I/P$ acts unipotently on $V=H^n(\overline X,\mathbb Q_\ell)$.
This allows to define a nilpotent monodromy operator $N=N_\ell$ on $V$ (basically the logarithm of a generator of $I/P$) satisfying
$$
N F\; =\; q F N
$$
for any (lift of a) {\it geometric} Frobenius $F$ in $G_K$.
Moreover one obtains an ascending \emph{monodromy filtration} $M_\cdot$ on $V$, characterized by the fact that $N M_i\subseteq M_{i-2}$ and that
\begin{equation}\label{eq.monodromy}
N^i\; :\; Gr^M_i V\; \stackrel\sim\longrightarrow \; Gr_{-i}^MV
\end{equation}
is an isomorphism for all $i\in \mathbb N$. More canonically, using the canonical $G_K/I$-isomorphism
\begin{equation}\label{eq.I/P}
I/P \cong \prod_{\ell\neq p}\,\zl(1)\,,
\end{equation}
$N$ can be seen as a $G_K$-equivariant map $V \rightarrow V(-1)$,
and \eqref{eq.monodromy} is a $G_K$-isomorphism
$$
Gr^M_i V\; \stackrel\sim\rightarrow \; Gr^M_{-i}V(-i)\,.
$$
By construction (and assumption that $P$ acts trivially), $Gr^M_iV$ is a $\mathbb Q_\ell$-representation of $G_K/I\cong G_k$.
There is the following conjecture (see \cite{De1} 8.1 and \cite{RZ}).

\medskip
\begin{conj}\label{conj.wt.monodr.l-adic} (\it monodromy weight conjecture)
The $G_k$-representation $Gr^M_iH^n(\overline{X},\ql)$ is pure of weight $n+i$.
\end{conj}

This was proved by Deligne in the equi-characteristic case $(\ch(K)>0)$, \cite{De6} (1.8.4),
at least if $X$ comes from a smooth projective scheme $\pi: \mathcal X \rightarrow U$
over an open subscheme $U$ of a curve $C$ over a finite field,
via passing to the henselization (or completion) at a closed point $x\in C-U$.
The general case was proved by T. Ito \cite{It}, by reducing it to Deligne's case.

\medskip
In the mixed characteristic case, i.e., the case of $p$-adic fields, very little is known on the
above conjecture. Using de Jong's (weak) resolution of singularities by alterations \cite{deJ},
and the result of M. Rapoport and Th. Zink \cite{RZ} for the case of semi-stable reduction,
one gets the conjecture for $\dim(X)\leq 2$, hence for $n \leq 2$ by a Lefschetz argument.
There is an analogous result in mixed Hodge theory, for a family of proper complex varieties
over a disk with degenerating fibre at 0, which produces the so-called limit mixed Hodge structure
on the generic fibre. The method in \cite{RZ} borrows from these techniques by using an \'etale
version of the vanishing cycles spectral sequence. The problem is that one needs a certain
non-degeneracy statement, which is known for the Hodge theory by some positivity, but is not known
for the $\ell$-adic case, except for surfaces where it is the classical Hodge index theorem.
It was proved by M. Saito \cite{SaM} that this non-degeneracy, and hence the monodromy weight conjecture,
would follow from Grothendieck's standard conjectures for varieties over finite fields.

\medskip
We discuss an implication for the fixed module. The following result is unconditional. Let
$\Gamma = G_K/I \cong G_k$ and $d = \dim(X)$.

\begin{thm}\label{thm.local.wts}
$H^n(\overline X,\ql)^I$ is a mixed $\Gamma$-module with weights in
$[\max(0,2(n-d)),\min(2n,2d)]$,
and $H^n(\overline X,\ql)_I$ is a mixed $\Gamma$-module with weights in $[\max(-2,2(n-d)),\min(2n,2d-2)]$.
In particular, $H^n(\overline X,\ql(r))^{G_K}=0$ for $r\notin [\max(0,n-d),\min(n,d)]$, and
$H^2(G_K,H^m(\overline X,\ql(s)))$ $= 0$ for $s \notin [\max(m+1-d,1),\min(d+1.m+1)]$.
\end{thm}

In fact, by de Jong's resolution of singularities \cite{deJ} one easily reduces to the case
where one has a regular proper model $\pi: \mathcal X \rightarrow \mathcal O_K$.
Then one has a long exact sequence
$$
\ldots \rightarrow H^n(\mathcal X_{ur},\ql) \rightarrow H^n(X_{ur},\ql) \rightarrow H^{n+1}_{\overline Y}(\mathcal X_{ur},\ql) \rightarrow \ldots,
$$
where $Y$ is the special fiber of $\pi$, $\overline Y = Y\times_k\overline{k}$ and the
subscript $ur$ denotes the base change to the maximal unramified extension of $K$, i.e.,
to the strict henselization of $\mathcal O_K$. By proper base change we have an isomorphism
$$
H^n(\mathcal X_{ur},\ql)\cong H^n(\overline Y,\ql)\,,
$$
and by Deligne's proof of the Weil conjecture \cite{De6} (3.3.8) the latter representation
is mixed with weights in $[\max(0,2(n-d)),2n]$. By the proof of cohomological
purity by O. Gabber and K. Fujiwara \cite{FuG} one has isomorphisms
$$
H^{n+1}_{\overline Y}(\mathcal X_{ur},\ql) \cong H_{2(d+1)-n-1}(\overline Y,\ql(d+1)) \cong H^{2d-n+1}(\overline Y,\ql(d+1))^\vee\,,
$$
where $V^\vee$ is the dual of a $\ql$-$G_k$-representation $V$, the group in the middle is $\ell$-adic \'etale
homology, and by Deligne (loc. cit.) the representation on the right is mixed with weights in $[2(n-d),\min(2d,2n)]$.
The first claim is now obvious from the exact sequence
$$
0 \longrightarrow H^1(I,H^{n-1}(\overline X,\ql)) \longrightarrow H^n(X_{ur},\ql)
\longrightarrow H^n(\overline X,\ql)^I \longrightarrow 0\,,
$$
which follows from the Hochschild-Serre spectral sequence and the fact that $I$ has $\ell$-coho\-mo\-lo\-gical
dimension 1. The second claim follows as well, using the isomorphisms
$$
H^1(I,H^{n-1}(\overline X,\ql)) \cong H^1(I/P,H^{n-1}(\overline X,\ql)) \cong H^{n-1}(\overline X,\ql))_I(-1)
$$
coming from the isomorphism \eqref{eq.I/P} and the fact that $P$ is a pro-p-group.
The third claim follows from the first claim, the isomorphism
\begin{equation}\label{eq.fix.modules}
H^n(\overline X,\ql)^{G_K} = (H^n(\overline X,\ql)^I)^\Gamma
\end{equation}
and Fact 2. The final claim is deduced by local duality.

\medskip
The monodromy weight conjecture gives a better bound.

\begin{cor}\label{cor.local.wts.l-adic}
If the monodromy weight conjecture holds, then

\smallskip
(a)\; $H^n(\overline X,\ql)^I$ is mixed with weights in $[\max(0,2n-2d),n]$ and
$H^n(\overline X,\ql)_I$ is mixed with weights in $[n,\min(2n,2d-2)]$.

\smallskip
(b) $H^n(\overline X,\ql(r))^{G_K} = 0$ for $r \notin [\max(0,n-d),\frac{n}{2}]$, and
$H^2(G_K,H^m(\overline X,\ql(s)) = 0$ for $s \notin [\frac m 2 + 1,\min(d+1,m+1)]$.
\end{cor}

In fact, by the construction of the monodromy filtration one has
$$
(H^n(\overline X,\ql)^I)^\Gamma \subset (M_0H^n(\overline X,\ql))^\Gamma\,,
$$
and the monodromy weight conjecture implies that $M_0H^n(\overline X,\ql)$ is a mixed $\Gamma$-module of
weights $\leq n$.
Similarly one has a surjection $H^n(\overline X,\ql)/M_{-1} \twoheadrightarrow H^n(\overline X,\ql)_I$
and monodromy weight conjecture implies that $H^n(\overline X,\ql)/M_{-1}$ is mixed of weights $\geq n$.
Thus (a) follows from Theorem \ref{thm.local.wts}. Now (b) follows as the last claim of \ref{thm.local.wts}.

\medskip
We end this section with the following speculation. In the equi-characteristic case Deligne proved
a more general result, where $C$ is a smooth curve over a finite field $k$, $U \subset C$ is a non-empty
open subscheme and $\ell$ is a prime invertible in $k$.

\medskip
\begin{thm}\label{thm.monodromy} (\cite{De6} (1.8.4)) Let $F$ be a smooth $\ql$-sheaf on $U$ which is
pure of weight $w$ (i.e., for each closed point $s\in U$ the stalk $F_{\overline{s}}$ is a pure $G_{k(s)}$-representation
of weight $w$). Let $x\in C-U$ be a closed point, let $K_x$ be the completion of the
function field of $C$ at $x$ (i.e., the fraction field of the completion of $\mathcal O_{C,x}$),
and let $\overline{\eta} = \Spec(\overline{K_x}) \rightarrow U$ be the corresponding geometric
point. If $M$ is the monodromy filtration on the $G_{K_x}$-representation $F_{\overline{\eta}}$
(the stalk of $F$ at $\overline{\eta}$), then $Gr^M_i(F_{\overline{\eta}})$ is pure of weight $w+i$.
\end{thm}

\begin{ques} Does the same hold if $U$ is an open subscheme of $\Spec(\mathcal O_K)$ for a number field $K$,
and $\ell\neq \ch(k(x))$?
\end{ques}

\vspace{0,5cm}

\section{Weights for $p$-adic cohomology over local fields}

\medskip
Let $K$ be a local field with finite residue field $k$ as before, but now assume that $\ch(K)=0$ and $\ell=p=\ch(k)$. Let again $X$ be a smooth projective variety over $K$. Even if $X$ has good reduction, the $\mathbb Q_p$-$G_K$-representation $V=H^n(\overline X,\mathbb Q_p)$ is not unramified, and there is no obvious way to see weights. However it was an insight of J.-M. Fontaine that one can associate a canonical object to $V$ on which one has a Frobenius and weights: Let $W = W(k)$ and $K_0 = Quot(W)$ as in section 3; $K_0$ is isomorphic to the maximal unramified extension of $\qp$ in $K$.
Fontaine \cite{Fo1} defined a certain ring $B_{cris}$ over $K_0$ and conjectured a comparison isomorphism
\begin{equation}\label{eq.comp.iso.crys}
H^n(\overline X,\mathbb Q_p)\otimes_{\mathbb Q_p}B_{cris}\; \cong\; H^n(\mathcal X_s/K_0)\otimes_{K_0} B_{cris}
\end{equation}
where $\mathcal X_s$ is the special fibre of a smooth projective model $\mathcal X/\Spec(\mathcal O_K)$ of $X$ over $\mathcal O_K$ as in section 4. Moreover, one has the following structures: $H^n(\overline X,\mathbb Q_p)$ is a $\mathbb Q_p$-$G_K$-representation and $H^n(\mathcal X_s/K_0)=H^n(\mathcal X_s/W)\otimes_WK_0$ is a so-called filtered $\varphi$-module over $K$, i.e., $H^n(\mathcal X_s/K_0)$ is a $K_0$-vector space with a Frobenius $\varphi$, and $H^n(\mathcal X_s/K_0)\otimes_{K_0}K$ has a descending filtration $F^\cdot$, via a canonical isomorphism
\begin{equation}\label{eq.comp.cry.dR}
H^n(\mathcal X_s/K_0)\otimes_{K_0}K\; \cong\; H^n_{dR}(X/K)
\end{equation}
and the Hodge filtration on the de Rham cohomology on the right. Now $B_{cris}$ has both structures -- it is a continuous $\mathbb Q_p$-$G_K$-representation and a $\varphi$-filtered module over $K$ -- and one can recover $H^n(\mathcal X_s/K_0)$ from $H^n(\overline X,\mathbb Q_p)$ and vice versa, viz., one has
\begin{equation}\label{eq.comp.1}
H^n(\mathcal X_s/K_0) \; =\; (H^n(\overline X,\mathbb Q_p)\otimes_{\mathbb Q_p}B_{cris})^{G_K}
\end{equation}
and
\begin{equation}\label{eq.comp.2}
H^n(\overline X,\mathbb Q_p) \; =\; (H^n(\mathcal X_s/K_0)\otimes_{K_0} B_{cris})^{\varphi=id}\cap F^0\, ,
\end{equation}
where the brackets have the mentioned three structures as well: $H^n(\overline X,\mathbb Q_p)\otimes_{\mathbb Q_p}B_{cris}$ has the diagonal $G_K$-action, and $H^n(\mathcal X/K_0)\otimes_{K_0}B_{cris}$ has the diagonal $\varphi$-action, and the usual tensor filtration after scalar extension to $K$ ($F^m(A\otimes B)=\sum_{i+j=m}F^iA\otimes F^jB$).

\medskip
In fact, Fontaine showed that there are mutual inverse category equivalences
\begin{equation}\label{eq.cat.equiv.crys}
\begin{array}{rcl}
\mbox{(crystalline $\mathbb Q_p$-$G_K$-representations)} & \rightarrow & \mbox{(admissible $\varphi$-filtered modules over $K$)}\\
V & \mapsto & D(V)=(V\otimes_{\mathbb Q_p}B_{cris})^{G_K}\\
(D\otimes_{K_0}B_{cris})^{\varphi=1} \cap F^0 & \leftmapsto & D
\end{array}
\end{equation}

where we refer to \cite{Fo1} and \cite{Pp} for a precise definition of both categories.
Fontaine's crystalline conjecture was shown by J.-M. Fontaine and W. Messing \cite{FM} for $\dim(X)<p$,
and by G. Faltings \cite{Fa1} in general.

\vspace{0,5cm}
We note that, in the geometric situation above, $D = H^n(\mathcal X_s/K_0)$ is pure of weight $n$, see section 3.

\vspace{0,5cm}
For $X/K$ with not necessarily good reduction, Fontaine and I arrived at the following conjecture
(\cite{Ja1} p. 347 and \cite{Fo2}), which was then subsequently proved by work of O. Hyodo and K. Kato
\cite{HyKa}, \cite{Ka2} and T. Tsuji \cite{Tsu}:
There is a ring $B_{st}$ over $K_0$ which is a continuous $\qp$-$G_K$-representation and also a
filtered $(\varphi,N)$-module over $K$, i.e.,
it has the same structures as $B_{cris}$, plus an operator $N$ such that
\begin{equation}\label{eq.N.phi}
N\varphi \;  = \;  p\varphi N\,.
\end{equation}
After possibly passing to a finite extension of $K$ (and over $K$ itself if $X$ has semi-stable reduction),
there is an isomorphism, compatible with $G_K,\varphi, F^\cdot$ and $N$,
\begin{equation}\label{eq.comp.iso.st}
H^n(\overline X,\mathbb Q_p)\otimes_{\mathbb Q_p}B_{st}\; \cong \; D^n\otimes_{K_0}B_{st}\, ,
\end{equation}
where $D^n$ is a finite-dimensional filtered $(N,\varphi)$-module over $K$. In fact, if $X$ has semi-stable reduction, then $D^n$ can be realized as the $n$-th log-crystalline cohomology of the special fibre $\mathcal X_s$ of a semistable model over $\mathcal O_K$.

\medskip
As in the $\ell$-adic case, the monodromy operator $N$ on $D^n$ -- which must be nilpotent by \eqref{eq.N.phi}
and finite-dimensionality of $D^n$ - allows to define a monodromy filtration $M_\cdot$ with
\begin{equation}\label{eq.monodromy.p-adic}
N^i\; : \; Gr^M_i D^n\; \stackrel\sim\rightarrow\; Gr_{-i}^MD^n(-i)\, ,
\end{equation}
and I conjectured (\cite{Ja1} p. 347) the following $p$-adic analogue of \eqref{conj.wt.monodr.l-adic}:

\begin{conj}\label{conj.wt.monodr.p-adic}({\it $p$-adic monodromy weight conjecture})
As a $\varphi$-module, $Gr_i^M D^n$ is pure of weight $n+i$ and has the same eigenvalues as
$Gr_i^MH^n(\overline X,\mathbb Q_\ell)$ for each $\ell\neq p$.
\end{conj}

The action of $\varphi$ is $\sigma$-semi-linear, but the $r$-th power of $\varphi$ (where $p^r=q$ is the cardinality
of the residue field of $K_0$ - and $K$) acts linearly, and hence there are well-defined
eigenvalues of $\varphi^r$, and the purity is defined via $q$ and these eigenvalues as in \eqref{eq.purity}.

\medskip
As in the $\ell$-adic mixed characteristic case, this conjecture is still wide open, except for
a few remarkable cases by T. Saito \cite{SaT1}, \cite{SaT2}.

\vspace{0,5cm}
There is an equivalence of categories
\begin{equation}\label{eq.cat.equiv.st}
(\mbox{semi-stable $\mathbb Q_p$-$G_K$-representations}) \quad \rightleftarrows \quad (\mbox{admissible filtered $(N,\varphi)$-modules over $K$})
\end{equation}
extending the one for the crystalline case \cite{Pp}.

\medskip
This is compatible with Tate twists on both sides. Therefore one has:

\begin{cor}\label{cor.local.wts.p-adic}
The monodromy weight conjecture \eqref{conj.wt.monodr.p-adic} would imply
$H^n(\overline X,\mathbb Q_p(r))^{G_K} = 0$ for $r \notin [\max(0,n-d),\frac{n}{2}]$, and
$H^2(G_K,H^m(\overline X,\qp(s)))=0$ for $s \notin [\frac m 2 + 1,\min(d+1,m+1)]$.
\end{cor}

In fact, the equivalence of categories would give an isomorphism
$$
\begin{array}{rcl}
H^n(\overline X,\mathbb Q_p(r))^{G_K} & \cong & Hom_{G_K}(\mathbb Q_p, H^n(\overline X,\mathbb Q_p(r)))\\
& \cong & Hom_{N,\varphi, F^\cdot}(K_0,D^n)\subseteq M_0D(r)^{\varphi=id}
\end{array}
$$
and $M_0D(r)$ would be mixed with weights in the interval $[\max(0,2n-2d)-2r,n-2r]$.

\medskip
Unconditionally one can use the so-called Hodge-Tate decomposition to show the
$p$-adic analogue of Theorem \ref{thm.local.wts}:

\begin{thm}\label{thm.local.wts.p-adic} (\cite{Ja1} p. 343 Cor. 5, or \cite{Sou} proof of Thm. 2 iii), plus
an easy improvement by hard Lefschetz)
One has $H^n(\overline X,\qp(r))^{G_K} = 0$ for $r\notin [\max(0,n-d),\min(n,d)]$.
\end{thm}

We end with a question similar to the one at the end of the previous section. Let $U$
be a non-empty open subscheme of $D=\Spec(\mathcal O_K)$ where $K$ is a number field,
let $p$ be a prime invertible on $U$, and let $F$ be a smooth $\qp$-sheaf on $U$
which is pure of weight $w$. Let $K_x$ be the completion of $K$ at a point $x \in D-U$
and assume that $p = \ch(k(x))$. Let $\overline \eta = \Spec(\overline K_x) \rightarrow U$
be the corresponding geometric point, and let $F_{\overline{\eta}}$ be the stalk of $F$
at $\overline \eta$, considered as a $\qp$-$G_{K_x}$-representation.

\begin{ques}
Is $F_{\overline{\eta}}$ a potentially semi-stable $\qp$-$K_{x}$-representation, i.e., a semi-stable
re\-pre\-sentation after restricting it to a finite field extension $K'$ of $K_{x}$? If $M_.$ is the
monodromy filtration on the filtered $(N,\phi)$-module $D$ associated to this $\qp$-$G_{K'}$-representation
via the category equivalence \eqref{eq.cat.equiv.st} over $K'$, is $Gr^M_iD$ pure of weight $w+i$
as a $\varphi$-module?
\end{ques}

\section{Weights and Galois cohomology}

\medskip
Let $K$ be a global field, i.e., a number field or a function field in one variable over a finite field.
Let $X$ be smooth projective variety of pure dimension $d$ over $K$. For many arithmetic applications the Galois cohomology groups
$$
H^1(K,H^n(\overline X,\mathbb Z_\ell(r)))
$$
need to be studied, together with the restriction map
\begin{equation}\label{eq.res.H1}
H^1(K, H^n(\overline X,\mathbb Z_\ell(r))) \; \rightarrow \prod\limits_v H^1(K_v, H^n(\overline X,\mathbb Z_\ell(r)))
\end{equation}
where $v$ runs over the places of $K$ and $K_v$ in the completion of $K$ at $v$, i.e., a local field. In fact, these data allow to define (generalized) Selmer groups which play a role in the conjectures of Birch, Swinnerton-Dyer, Beilinson, Bloch and Kato on generalized class number formulae and special values of $L$-functions.

\medskip
Therefore it is interesting to study the kernel of \eqref{eq.res.H1}. By Poitou-Tate duality, and Poincar\'e duality
$$
H^n(\overline X,\mathbb Z_\ell(r))\times H^{2d-n}(\overline X,\mathbb Q_\ell/\mathbb Z_\ell(d-r))\; \rightarrow\; H^{2d}(\overline X,\mathbb Q_\ell/\mathbb Z_\ell(d))\; \stackrel{Tr}{\rightarrow} \mathbb Q_\ell/\mathbb Z_\ell\,,
$$
the kernel of \eqref{eq.res.H1} is dual to the kernel of
\begin{equation}\label{eq.res.H2}
H^2(K,H^{m}(\overline X,\mathbb Q_\ell/\mathbb Z_\ell(s)))\; \rightarrow\; \mathop\oplus\limits_vH^2(K_v, H^{m}(\overline X,\mathbb Q_\ell/\mathbb Z_\ell(s)))
\end{equation}
for $m=2d-n$ and $s=d-r+1$.
It turns out that the kernel can be controlled by weights. In fact, one has the following cohomological Hasse principle:

\begin{thm}\label{thm.coh.Hasse.pr} (\cite{Ja4} Thm. 1.5)  Let $A$ be a discrete $G_K$-module which is a cofinitely generated divisible torsion $\mathbb Z_\ell$-module (i.e., isomorphic to $(\mathbb Q_\ell/\mathbb Z_\ell)^r$ for some $r$, as a $\mathbb Z_\ell$-module). If $A$ is mixed of weights $\neq -2$, the restriction map
$$
H^2(K,A)\; \stackrel\sim\rightarrow\; \mathop\oplus\limits_v H^2(K_v,A)
$$
is an isomorphism. Here $A$ is called mixed of weights $w_1,\ldots , w_m$, if this holds for the $r$-dimensional $\mathbb Q_\ell$-representation $V_\ell(A)=T_\ell(A)\otimes_{\mathbb Z_\ell}\mathbb Q_\ell$, where $T_\ell(A)=\lim\limits_{\leftarrow i}A[\ell^i]$ is the inverse
limit over the modules $A[\ell^i] = \{a\in A \mid \ell^ia=0\}$ via the $\ell$-multiplications.
\end{thm}

It follows that \eqref{eq.res.H2} has finite kernel for $m-2s\neq -2$, i.e., $m\neq 2(s+1)$ because $m-2s$ is the (pure) weight of
$DivH^{m}(\overline X,\mathbb Q_\ell/\mathbb Z_\ell(s))$, where $DivA$ is the maximal divisible module
of a torsion Galois module $A$.

\medskip
By Poitou-Tate duality, the kernel of \eqref{eq.res.H1} is finite for $n\neq 2r$.

\medskip
Combined with the local vanishing results from sections 4 and 5 we obtain

\begin{thm}\label{thm.vanish.H2}
One has $H^2(K, DivH^m(\overline X,\qzl(s)))=0$ for $s\notin [\max(-d-1,1),\min(d+1,m+1)]$.
\end{thm}

In fact, the local groups $H^2(K_v,DivH^m(\overline X,\qzl(s)))$ are a quotient of
$H^2(K_v,H^m(\overline X,\ql(s)))$, which is dual to $H^n(\overline X,\ql(r))^{G_K}$
for $n=2d-m$ and $r=d-s+1$. Now we can apply Theorems \ref{thm.local.wts} and
\ref{thm.local.wts.p-adic}.

\medskip
By the same argument, the local monodromy weight conjectures would imply:

\begin{conj}\label{conj.vanish.H2}
$H^2(K,Div (H^m(\overline X,\qzl(s)))=0$ for $s\notin [\frac m 2+1,\min(d+1,m+1)]$.
\end{conj}

\medskip
Arithmetic applications often need the Galois cohomology groups
$$
H^i(G_S, H^n(\overline X,\mathbb Z_\ell(r)))\, .
$$
Here $S$ is a finite set of places $v$ of $K$ and $G_S=Gal(k_S/k)$ is the Galois group of the maximal extension $k_S$ of $k$ which is unramified outside $S$. Moreover $S$ should contain all archimedean places and all places $v|\ell$, and should be large enough so that the action of $G_K$ on $H^n(\overline X,\mathbb Z_\ell(r))$ factors through $G_K\twoheadrightarrow G_S$, i.e., such that $H^n(\overline X,\mathbb Z_\ell(r))$ is unramified outside $S$. Such an $S$ always exists. In fact, there in open subscheme $U\subseteq C=\Spec(\mathcal O_K)$ (if $K$ is a number field) or $U\subseteq C$ (if $K$ is the function field of a smooth projective curve $C$ over a finite field $\mathbb F_q$), such that $X/K$ has a smooth projective model $\pi:\mathcal X\rightarrow U$. As we have seen in section 2, $H^n(\overline X,\mathbb Z_\ell(r))$ is then unramified at all places $v$ corresponding to points $x\in U$ with $v\nmid \ell$, so we may take $S$ as the set of places not corresponding to these $x$, which are only finitely many. Then, after possibly shrinking $U$, we can write the above Galois cohomology group also as \'{e}tale cohomology
$$
H^i(U, R^n\pi_\ast\mathbb Z_\ell(r))\, .
$$
There is the following vanishing conjecture (\cite{Ja1} Conj. 1).

\begin{conj}\label{conj.Conj1} Let $K$ be a number field. Then one has
$$
H^2(G_S,H^m(\overline X,\mathbb Q_\ell(s)) = 0 \quad\mbox{  for  }\quad s \notin [\frac{m}{2}+1,min(d+1,m+1)]\,.
$$
\end{conj}

The statement in loc. cit. uses the condition $s\notin [\frac m 2 + 1,m+1]$; the refinement comes
from a simple Lefschetz argument.

\medskip
Obviously, Conjecture \ref{conj.Conj1} would follow from the monodromy weight conjectures and the
following one:

\begin{conj}\label{conj.global.local.GS}
If $K$ is a number field, then the map
$$
H^2(G_S,H^m(\overline X,\ql(s))) \longrightarrow \oplus_{v\in S}\; H^2(K_v,H^m(\overline X,\ql(s)))
$$
is injective for $m-2s\neq -1$.
\end{conj}

\medskip
By Poitou-Tate duality, this is equivalent to the injectivity of
\begin{equation}\label{eq.global.local.GS.H1}
H^1(G_S,H^n(\overline X,\ql(r))) \longrightarrow \oplus_{v\in S}\; H^1(K_v,H^n(\overline X,\ql(r)))
\end{equation}
for $n-2r\neq -1$.
The two conjectures above seem rather difficult, but they are partly motivated by the fact that
the same statements are true over global function fields. As for \ref{conj.Conj1}
one has the following more general result:

\begin{thm} \label{thm.vanishing.U} (\cite{Ja1} p. 335 Thm. 2) Let $U$ be a smooth curve over a finite field $\mathbb F_q$, let $\ell$ be prime, $\ell\neq \ch(\mathbb F_q)$, and let $F$ be a smooth (=twisted constant) $\mathbb Q_\ell$-sheaf on $U$. Assume that $F$ is mixed of weights $\geq 0$
(has a filtration with pure quotients of weights $\geq 0$) or that $F^\vee(2)$ is entire (\cite{De6} 1.22), where $F^\vee$ is the dual of $F$. Then $H^2(U,F)=0$.
\end{thm}

\medskip
 The statement in Conjecture \ref{conj.Conj1} is a special case, because $H^n(\overline X,\mathbb Q_\ell(j))$ regarded as the smooth $\mathbb Q_\ell$-sheaf $R^n\pi_\ast\mathbb Q_\ell(j)$ for $\pi:\mathcal X\rightarrow U$ as above, has weights $\geq 0$ if $n-2j\geq 0$, i.e., $n+1>2j$, and $H^n(\overline X,\mathbb Q_\ell(j))^\vee(2)\cong H^{2d-n}(\overline X,\mathbb Q_\ell(d-j))(2)\cong H^n(\overline X,\mathbb Q_\ell(n-j+2))$ (Poincar\'{e} duality and hard Lefschetz) is entire for $n-j+2\leq 0$, i.e., $n+1<j$.

\medskip
The statement in \ref{conj.global.local.GS} for global function fields follows from the more general:

\medskip
\begin{thm}\label{thm.global.local.U} Let $U$ be a smooth curve over a finite field $\mathbb F_q$,
let $\ell$ be prime, $\ell\neq \ch(\mathbb F_q)$, and
let $F$ be smooth $\mathbb Q_\ell$-sheaf of weight $w$ on $U$. Then the restriction map
$$
H^2(U,F)\; \rightarrow\; \mathop\oplus\limits_{v\in S}H^2(K_v,F)
$$
is injective for $w\neq -1$.
\end{thm}

\textbf{Proof }
Let $j:U\hookrightarrow C$ be an open immersion into a smooth projective curve $C$.
Then the above map can be identified with the map $\alpha$ in the long exact cohomology sequence
\begin{equation}\label{eq.rel.coh.seq}
\ldots \rightarrow H^2(C,j_!F)\; \stackrel\beta\rightarrow \; H^2(U,F)\; \stackrel\alpha\rightarrow\; \mathop\oplus\limits_{v\in S=C-U} H^3_v(C,j_!F)\; \rightarrow\; H^3(C,j_!F)\rightarrow\ldots
\end{equation}
and the map $\beta$ factorizes through $H^2(C,j_\ast F)$ which sits in an exact sequence
\begin{equation}\label{eq.HS.curves}
0\; \rightarrow\; H^1(\overline C, j_\ast F)_\Gamma\; \rightarrow \; H^2(C, j_\ast F)\; \rightarrow \; H^2(\overline C,j_\ast F)^\Gamma\;\rightarrow \; 0\, .
\end{equation}
Now Deligne has proved that $H^i(\overline C, j_\ast F)$ is pure of weight $w+i$, so $H^2(C, j_\ast F)$ vanishes if $w+1\neq 0$ and $w+2\neq 0$. The case $w=-2$ follows with some extra argument (loc. cit. and \cite{Ra} Thm. 4.1).

\bigskip
We dare to state the following conjecture which would contain conjectures \ref{conj.Conj1} and
\ref{conj.global.local.GS}.

\begin{conj}\label{conj.vanishing.U} Theorems \ref{thm.vanishing.U} and \ref{thm.global.local.U}
also holds for an open subscheme $U \subset \Spec(\mathcal O_K)$, where $K$ is a number field and
$\ell$ is invertible on $U$.
\end{conj}

\medskip
The exact sequence \eqref{eq.rel.coh.seq} and the factorization
$$
\beta\; :\; H^2(C, j_! F)\; \rightarrow\; H^2(C, j_\ast F)\; \rightarrow\; H^2(U,F)
$$
also exist for a number field $K$ and $C=\Spec(\mathcal O_K)$. Thus we are tempted to state

\begin{conj}\label{conj.j!} Let $\ell$ be a prime and let $U\subseteq C$ be an open subscheme such that $\ell$ is invertible on $U$. Let $F$ be a smooth $\mathbb Q_\ell$-sheaf on $U$ which is mixed of weights $\neq -1$. Then
$$
H^2(C, j_\ast F)=0\, .
$$
\end{conj}

The problem is that there are no obvious analogues of the groups $H^i(\overline C, j_! F), H^i(\overline C, j_\ast F)$ etc. in the number field case.
The most common analogue of the morphism $\overline C\rightarrow C$, a pro-\'{e}tale covering with Galois group $\Gamma\cong \hat{\mathbb Z}$, would be the $\ell$-cyclotomic extension $\overline U\rightarrow U$ with Galois group $\mathbb Z^\times_\ell$. One even can consider the corresponding covering $\overline C\rightarrow C$, but it is not an \'{e}tale covering (so we do not have a Hochschild-Serre spectral sequence as in \eqref{eq.HS.curves} above), the associated cohomology groups do not have finiteness properties, and there does not seem to be a theory of weights on them.

\vspace{0,5cm}
But we should remember that the situation was the same at the local places $v$ with $v|\ell$. There the weights only became visible after applying Fontaine's comparison functors. This leads to the following.

\begin{ques}
Does there exist a global analogue of Fontaine's functors over number fields?
\end{ques}

We recall that, although Iwasawa theory also exists and is useful over $p$-adic fields, even there the obtained modules are only finitely generated over the Iwasawa algebra, and do not show weights. Conversely, Fontaine's theory of $B_{cris}$ and $B_{st}$ does not see the cyclotomic extension, at least not directly. Nevertheless a certain link between Fontaine's theory and $\mathbb Z_p$-extensions is given by the theory of $(\varphi,\Gamma)$-modules and the fields of norms.

\begin{ques}
Do there exist global analogues of these?
\end{ques}

\vspace{0,5cm}
\section{Application: Hasse principles for function fields}

\medskip
The cohomological Hasse principle in the previous section (Theorem \ref{thm.coh.Hasse.pr})
led to a proof of the following Hasse principle conjectured by K. Kato \cite{Ka1}, and proved by him for $d=1$.

\medskip
Let $K$ be a global field, let $F/K$ be a function field in $d$ variables which is primary (i.e., such that $K$ is separably closed in $F$), and let $\ell$ be a prime. For every place $v$ of $K$ let $F_v$ be the corresponding function field over $K_v$: If $F=K(V)$ for a geometrically integral variety $V$ over $K$, then $F_v=K_v(V_v)$, where $V_v=V\times_K K_v$.

\begin{thm}\label{thm.Hasse.ff} (\cite{Ja4} Thm. 2.7) Let $K$ be a number field. Then the map
$$
H^{d+2}(F,\mathbb Q_\ell/\mathbb Z_\ell(d+1))\; \rightarrow \; \mathop\oplus\limits_v H^{d+2}(F_v,\mathbb Q_\ell/\mathbb Z_\ell(d+1))
$$
is injective.
\end{thm}

As in the classical case $d=0$ (which corresponds to the classical Hasse principle for the Brauer group of $K$),
or the case $d=1$ (see \cite{Ka1}, appendix by Colliot-Th\'el\`ene),
this Hasse principle has applications to quadratic forms. In fact, it implies that the Pythagoras
number of $F$ is bounded by $2^{d+1}$ if $d\geq 2$ (\cite{CTJ}).

\medskip
For the proof of Theorem \ref{thm.Hasse.ff} one first shows that, via the Hochschild-Serre spectral
sequence, the kernel of the map above is isomorphic to the kernel of the map
$$
H^2(K, H^d(F\overline K,\mathbb Q_\ell/\mathbb Z_\ell(d+1)))\;\rightarrow\; \mathop\oplus\limits_v H^2(K_v, H^d(F\overline K,\mathbb Q_\ell/\mathbb Z_\ell(d+1)))
$$
where $F\overline K=\overline K(\overline V)$, with $\overline V=V\times_K \overline K$, is the corresponding function field over $\overline K$.

\medskip
Now we have
$$
H^d(F\overline K,\mathbb Q_\ell/\mathbb Z_\ell(d+1)) \; =\; \lim\limits_\rightarrow H^d(\overline U,\mathbb Q_\ell/\mathbb Z_\ell(d+1))\, ,
$$
where $U\subseteq V$ runs over the smooth open subvarieties. Thus it suffices to show the injectivity of
$$
H^2(K, H^d(\overline U,\mathbb Q_\ell/\mathbb Z_\ell)(d+1))\; \rightarrow\; \mathop\oplus\limits_v H^2(K_v, H^d(\overline U,\mathbb Q_\ell/\mathbb Z_\ell)(d+1))
$$
for all affine open smooth $U\subseteq V$, or at least a cofinal set of them. To apply Theorem \ref{thm.coh.Hasse.pr}, we investigate the weights of $N=H^d(\overline U,\mathbb Q_\ell/\mathbb Z_\ell)(d+1)$. Now $H^d(\overline U,\mathbb Q_\ell/\mathbb Z_\ell)$ is divisible (by weak Lefschetz, see
\eqref{eq.weak.Lefsch} below) and mixed of weights $d,\ldots , 2d$, so $N$ is mixed of weights $-2d-2,-2d-1,\ldots, -3,-2$. We cannot apply Theorem \ref{thm.coh.Hasse.pr} directly, because the weight $-2$ occurs. But one can use the weight filtration to show that it suffices to consider the weight $-2$ quotient, which sits at the top of the ascending weight filtration.

\medskip
There are many weight $-2$ modules for which the Hasse principle fails -- basically this amounts to tori over $K$ for which the analogous Hasse principle is wrong. So we have to study $Gr^W_{-2}N$ carefully.

\vspace{0,5cm}
For this one uses resolution of singularities, which holds over fields of characteristic zero by the work of Hironaka: There is a good compactification $U\subset X\supset Y=X-U= \textstyle\bigcup^N_{i=1}Y_i$ as in section 1 (but over $K$), and we use the associated weight spectral sequence. But in general it has many terms and many non-vanishing differentials.

\vspace{0,5cm}
The next observation is that one can greatly simplify the situation by the weak Lefschetz theorem. It says that
\begin{equation}\label{eq.weak.Lefsch}
H^a(\overline U,\mathbb Z/\ell^n\mathbb Z)=0\quad\mbox{for}\quad a >d
\end{equation}
if $U$ is affine of dimension $d$. Now one chooses $U$ affine and a smooth hyperplane section $Y_{N+1}\subset X$ which intersects the normal crossing divisor $Y$ transversally. Then $Z = Y \cup Y_{N+1}$ is a divisor with normal crossings, $X^0=X-Y_{N+1}$, $U^0=U-Y_{N+1}=X-Z$ and $Y^0=Y-Y_{N+1}$ are affine, many terms in the spectral sequence vanish and one gets an exact sequence
\begin{equation}\label{eq.Gr.W-2}
0\;\rightarrow\; Gr^W_{-2}H^d(\overline{U^0},\mathbb Q_\ell/\mathbb Z_\ell(d+1))\;\rightarrow \;H^0(\overline{Z^{[d]}},\mathbb Q_\ell/\mathbb Z_\ell(1))\;\stackrel\gamma\rightarrow\; H^2(\overline{Z^{[d-1]}},\mathbb Q_\ell/\mathbb Z_\ell(2))\, .
\end{equation}
Here $H^0(-)$ and $H^2(-)$ are induced Galois modules, and a finer analysis shows that one has a Hasse principle for $Gr^W_{-2}H^d(\overline{U^0},\mathbb Q_\ell/\mathbb Z_\ell(d+1))$.
Since one can choose a cofinal set of these $U^0$, one obtains Theorem \ref{thm.Hasse.ff}.

\vspace{0,5cm}
Theorem \ref{thm.Hasse.ff} extends to global fields $K$ of positive characteristic if one assumes the
existence of good compactifications over the perfect hull of $K$. Due to recent results on resolution
of singularities \cite{CP1}, \cite{CP2}, \cite{CJS}, this holds for $d\leq 3$.

\vspace{0,5cm}
\section{Application: Hasse principles for smooth projective varieties over global fields}

\medskip
Now we consider the cokernel of the map in Theorem \ref{thm.Hasse.ff}. For any variety $X$ over
a global or local field, Kato defined a complex of Bloch-Ogus type $C^{2,1}(X,\mathbb Q_\ell/\mathbb Z_\ell)$:
\begin{equation}\label{eq.Kato.cx}
\begin{array}{ccccl}
\ldots&  \rightarrow& \mathop\oplus\limits_{x\in X_a} H^{a+2}(k(x), \mathbb Q_\ell/\mathbb Z_\ell(a+1))& \rightarrow&\mathop\oplus\limits_{x\in X_{a-1}} H^{a+1}(k(x), \mathbb Q_\ell/\mathbb Z_\ell(a))\rightarrow\ldots\\\\
\ldots & \rightarrow & \mathop\oplus\limits_{x\in X_1} H^3(k(x),\qzl(2)) & \rightarrow &\; \mathop\oplus\limits_{x\in X_0} \; H^2(k(x),\qzl(1))\, .
\end{array}
\end{equation}
Here $X_a$ is the set of points $x\in X$ of dimension $a$, $k(x)$ is the residue field of $x$, and the term with $X_a$ is placed in (homological) degree $a$. Then we have

\begin{thm}\label{thm.Hasse.X} (\cite{Ja4} Thm. 4.8)
Let $K$ be a number field, let $X$ be a connected smooth projective variety over $K$, and, for any place $v$ of $K$, let $X_v=X\times_K K_v$ be the corresponding variety over $K_v$. Then the restriction map
$$
C^{2,1}(X, \mathbb Q_\ell/\mathbb Z_\ell)\; \rightarrow \; \mathop\oplus\limits_v C^{2,1}(X_v,\mathbb Q_\ell/\mathbb Z_\ell)
$$
is injective, and for the cokernel $C'(X,\mathbb Q_\ell/\mathbb Z_\ell)$ one has
\begin{equation}\label{eq.Hasse.X}
H_a(C'(X,\mathbb Q_\ell/\mathbb Z_\ell))=\left\{\begin{array}{ccc} 0 & , & a>0,\\\mathbb Q_\ell/\mathbb Z_\ell& , & a=0\end{array}\right.\, .
\end{equation}
\end{thm}

This was conjectured by Kato (\cite{Ka1} Conj. 0.4) (for arbitrary global fields), is the classical sequence of Brauer groups
\begin{equation}\label{eq.Brauer}
0\rightarrow Br(K)\rightarrow \mathop\oplus\limits_v Br(K_v)\rightarrow \mathbb Q_\ell/\mathbb Z_\ell\rightarrow 0
\end{equation}
for $d=0$ and $X=\Spec(K)$, and was proved by Kato for $d=1$ in \cite{Ka1}. Theorem \ref{thm.Hasse.X} extends to a
global function field $K$, if resolution of singularities holds over the perfect hull $K'$ of $K$.

\bigskip
The first claim in Theorem \ref{thm.Hasse.X} easily follows from Theorem \ref{thm.Hasse.ff} above,
because the components of $C^{2,1}$ involve exactly the Galois cohomology groups considered in
\ref{thm.Hasse.ff} for all residue fields of $X$.

\medskip
For the proof of the second claim (and in fact even for the proof of Theorem \ref{thm.Hasse.ff})
it is useful to consider the henselizations $K_{(v)}$ of $K$ rather than the completions, and
the factorization
\begin{equation}\label{eq.hens.compl}
C^{2,1}(X,\qzl) \mathop{\longrightarrow}\limits^{\alpha} \mathop\oplus\limits_v \;C^{2,1}(X_{(v)},\qzl) \longrightarrow \mathop\oplus\limits_v\;C^{2,1}(X_v,\qzl)\,,
\end{equation}
where $X_{(v)} = X\times_KK_{(v)}$. This is possible because of the following rigidity result.

\begin{thm}\label{thm.rigidity}
The map of complexes $C^{2,1}(X_{(v)},\qzl) \longrightarrow C^{2,1}(X_v,\qzl)$ is a quasi-isomorphism
(i.e., induces an isomorphism in the homology) for all $v$.
\end{thm}

By this property one may replace the complex $C'(X,\qzl)$ with the complex $\overline C(X,\qzl)$
which is the cokernel of the map $\alpha$ in \eqref{eq.hens.compl}, and show \eqref{eq.Hasse.X}
for this complex.

\bigskip
Note that the complexes $\Cb(X,\mathbb Q_\ell/\mathbb Z_\ell)$ exist for arbitrary (not necessarily
smooth projective) varieties. Moreover, like for the complexes $C^{2,1}(X,\mathbb Q_\ell/\mathbb Z_\ell)$
and $C^{2,1}(X_{(v)},\mathbb Q_\ell/\mathbb Z_\ell)$, one has canonical short exact sequences of complexes
$$
0\rightarrow \Cb(Y,\mathbb Q_\ell/\mathbb Z_\ell)\rightarrow \Cb(X,\mathbb Q_\ell/\mathbb Z_\ell)\rightarrow \Cb(V,\mathbb Q_\ell/\mathbb Z_\ell)\rightarrow 0
$$
for $Y\subseteq X$ closed with open complement $V=X-Y$, because $X_a$ is the disjoint union of $Y_a$ and $V_a$ for all $a$.

\medskip
This gives rise to a so-called (Borel-Moore type) homology theory
$$
X\;  \mapsto\; \Hb_a(X,\mathbb Q_\ell/\mathbb Z_\ell):=H_a(\Cb(X,\mathbb Q_\ell/\mathbb Z_\ell))
$$
on the category $V_K^\ast$ of all varieties over $K$, with proper morphisms as morphisms, i.e., a sequence of covariant functors (for proper morphisms) from $V^\ast_K$ to abelian groups together with long exact sequences
$$
\ldots\rightarrow \Hb_a(Y,\mathbb Q_\ell/\mathbb Z_\ell)\stackrel{i_\ast}\rightarrow \Hb_a(X,\mathbb Q_\ell/\mathbb Z_\ell)\stackrel{j^\ast}\rightarrow \Hb_a(V,\mathbb Q_\ell/\mathbb Z_\ell)\stackrel\delta\rightarrow \Hb_{a-1}(Y,\mathbb Q_\ell/\mathbb Z_\ell)\rightarrow\ldots
$$
for every closed immersion $i: Y\hookrightarrow X$ with open complements $j:V\hookrightarrow X$, such that these long exact sequence are compatible with proper maps and the additional morphisms $j^\ast$ for open immersions, in an obvious way.

\vspace{0,5cm}
Next we observe that we can compute $\Hb_a$ rather well if we have resolution of singularities and if we \emph{assume} that property \eqref{eq.Hasse.X} holds for smooth projective varieties. For example, if $U$ is a connected smooth quasi-projective variety over $K$ and
\begin{equation}\label{eq.good.compact}
U\subset X\supset Y=\textstyle\bigcup^N_{i=1} Y_i
\end{equation}
is a good compactification, then one has an analogue of the weight spectral sequence
which becomes now very simple because of \eqref{eq.Hasse.X}. In fact it would follow that $\Hb_a(U,\mathbb Q_\ell/\mathbb Z_\ell)$ is the $a$-th homology of the complex
\begin{equation}\label{eq.wt.homology}
\mathbb Q_\ell/\mathbb Z_\ell^{\pi_0(Y^{[d]})}\rightarrow \mathbb Q_\ell/\mathbb Z_\ell^{\pi_0(Y^{[d-1]})}\rightarrow \hspace{1cm}\ldots \rightarrow \mathbb Q_\ell/\mathbb Z_\ell^{\pi_0(Y^{[1]})}\rightarrow \mathbb Q_\ell/\mathbb Z_\ell
\end{equation}
where $d=\dim(U)$ and the differentials have the obvious combinatoric description.

\vspace{0,5cm}
The proof of \eqref{eq.Hasse.X} then is based on the following result.

\begin{thm}\label{thm.homology.iso}
(i) There is a homology theory $H^W_\cdot(-,\mathbb Q_\ell/\mathbb Z_\ell)$ on $V^\ast_K$ which
has the property that for $U\subset X\supset Y$ as in \eqref{eq.good.compact},
$H_a^W(U, \mathbb Q_\ell/\mathbb Z_\ell)$ is the $a$-th homology of \eqref{eq.wt.homology}.

\medskip
(ii) There is a morphism of homology theories
$$
\varphi:\;\Hb_\cdot(-,\mathbb Q_\ell/\mathbb Z_\ell)\rightarrow H_\cdot ^W(-,\mathbb Q_\ell/\mathbb Z_\ell)\, .
$$

\medskip
(iii) This morphism is an isomorphism.
\end{thm}

Obviously this implies \eqref{eq.Hasse.X} for $\overline{C}$, because this condition holds for the homology $H^W(-,\qzl)$
by \ref{thm.homology.iso} (i), hence for $\Hb(-,\mathbb Q_\ell/\mathbb Z_\ell)$ by 8.3 (iii).

\medskip
The construction of (i) and (ii) depends on resolution of singularities and will be discussed
in the next section.

\medskip
The crucial point of (iii) is that obviously one can show it by induction on dimension and localization, i.e.,
by showing that for any integral variety $V$ with function field $K(V)$ the morphism
\begin{equation}\label{eq.varphi.generic}
\Hb_a(K(V),\mathbb Q_\ell/\mathbb Z_\ell)\; \rightarrow \; H^W_aÊ(K(V), \mathbb Q_\ell/\mathbb Z_\ell)
\end{equation}
is an isomorphism for all $a$, where
\begin{equation}\label{eq.lim.V'}
\Hb_a(K(V), \mathbb Q_\ell/\mathbb Z_\ell)\; =\; \lim\limits_{\rightarrow\,V'}\, \Hb_a(V',\mathbb Q_\ell/\mathbb Z_\ell)
\end{equation}
(with the limit running over all open subvarieties $V'\subset V$), similarly for $H^W$.

\medskip
Concerning \eqref{eq.varphi.generic}, we note that by definition of $\Cb(-,\qzl)$ one has
\begin{equation}\label{eq.vanishing.Hb}
\Hb_a(K(V),\qzl) = 0 \quad\quad \mbox{ for } a \neq d = \dim(V)\,.
\end{equation}
and
\begin{equation}\label{eq.Hb0}
\begin{array}{rcl}
&& \Hb_d(K(V),\qzl) =  \Cb_d(K(V),\qzl) \\\\
& = & \coker[H^{d+2}(K(V),\qzl(d+1)) \rightarrow \mathop\oplus\limits_v\,H^{d+2}(K(V)_{(v)},\qzl(d+1))] \\\\
& = & H^d(\overline{K}(\overline{V}),\qzl(d))_{G_K}\,,
\end{array}
\end{equation}
where $\overline{V}=V\times_K\overline{K}$ and $A_{G_K}$ is the cofixed module (maximal quotient with
trivial action) of a discrete $G_K$-module $A$. The last isomorphism follows easily from Poitou-Tate duality
(Here one uses that we consider henselizations $K_{(v)}$ instead of completions $K_v$).

\medskip
It remains to show the same properties for the weight homology $H^W_{\ldot}(-,\qzl)$.
Property \eqref{eq.vanishing.Hb} is shown by a Lefschetz argument:
If $U \subset X \supset Y$ is as above, with $U$ of pure dimension $d$,
then by Bertini's theorem we can find a smooth hyperplane section $Y_{N+1}$ of $X$
intersecting $Y$ transversally and having the property that
$$
\pi_0(Y^{[i]}\cap Y_{N+1}) \mathop{\longrightarrow}\limits^{\sim} \pi_0(Y^{[i]})
$$
is an isomorphism for all $i \leq d-2$ (viz., where $\dim(Y^{[i]}) \geq 2$) and that
$$
\pi_0(Y^{[d-1]}\cap Y_{N+1}) \twoheadrightarrow \pi_0(Y^{[d-1]})
$$
is a surjection (note $\dim(Y^{[d-1]}) = 1$). Let $U^0 = U-U\cap Y_{N+1}$, and note that
$U^0 = X - Z$, where $Z := Y\cup Y_{N+1}$ is again a simple normal
crossings divisor on $X$, by the transversality of $Y_{N+1}$ and $Y$.

Thus the commutative diagram
\begin{equation}\label{eq.double.cpx}
\begin{matrix}
 &  & (\qzl)^{\pi_0(Y^{[d-1]}\cap Y_{N+1})} & \longrightarrow & (\qzl)^{\pi_0(Y^{[d-2]}\cap Y_{N+1})} & \rightarrow & \ldots \\
 &  &  \downarrow                         &                 & \downarrow\rlap{$\cong$}  & & \downarrow\rlap{$\cong$} \\
(\qzl)^{\pi_0(Y^{[d]})} & \rightarrow & (\qzl)^{\pi_0(Y^{[d-1]})} & \rightarrow & (\qzl)^{\pi_0(Y^{[d-2]})} & \rightarrow & \ldots \,,
\end{matrix}
\end{equation}
in which the first vertical arrow is a surjection, shows that $H_a(U^0,\qzl) = 0$ for $a\neq d$,
because this is the $a$-th homology of the total complex associated to the double complex \eqref{eq.double.cpx}.
Since these $U^0$ are cofinal in the inductive limit \eqref{eq.lim.V'}, condition
\eqref{eq.vanishing.Hb} follows for $H^W_.$.
As for property \eqref{eq.Hb0} one has
$$
H^W_d(U^0,\qzl)
%& = & \ker((\qzl)^{\pi_0(Y^{[d]})} \oplus (\qzl)^{\pi_0(Y^{[d-1]}\cap Y_{N+1})} \longrightarrow
%(\qzl)^{\pi_0(Y^{[d-1]})} \oplus (\qzl)^{\pi_0(Y^{[d-2]}\cap Y_{N+1})}) \\
= \ker((\qzl)^{\pi_0(Z^{[d]})} \mathop{\longrightarrow}\limits^{\alpha^W} (\qzl)^{\pi_0(Z^{[d]})})
$$
by Theorem \ref{thm.homology.iso} (1). Now one shows the bijectivity of
$ \varphi: \Hb_d(U^0,\qzl) \mathop{\longrightarrow}\limits^{\sim} H^W_d(U^0,\qzl)$
by showing:

\begin{prop} There is a commutative diagram with exact rows
$$
\begin{CD}
0 @>>> H^d(\overline{U^0},\qzl(d))_{G_K} @>>> H^0(\overline{Z^{[d]}},\qzl)_{G_K} @>>> H^2(\overline{Z^{[d-1]}},\qzl)_{G_K} \\
@. @VV{\varphi}V @V{\cong}V{\varphi}V @ V{\cong}V{\varphi}V \\
0 @>>> H^W_d(U^0,\qzl) @>>> (\qzl)^{\pi_0(Z^{[d]})} @>>> (\qzl)^{\pi_0(Z^{[d-1]})}.
\end{CD}
$$
\end{prop}

\medskip
The upper row comes from the weight spectral sequence for $H^\ast(\overline{U^0},\qzl)$
(compare \eqref{eq.Gr.W-2}).

\section{Hypercoverings, hyperenvelopes, and weight complexes}

Hypercoverings were used by Deligne to define weight filtrations and mixed Hodge
structures on the cohomology of arbitrary complex algebraic varieties $X$. A covering of $X$ is
a surjective proper morphism $X' \rightarrow X$. It is called a smooth covering if $X'$ is smooth.
A simplicial variety $X.$ over $X$ (which is a simplicial object in the category of varieties
over $X$, or, equivalently, a morphism of simplicial varieties $X. \rightarrow X$
where $X$ also stands for the constant simplicial variety associated to $X$) is called a hypercovering
if for all $n\geq 0$, the morphism
\begin{equation}\label{eq.cosk}
X_n \longrightarrow (cosk_{n-1}^Xsk_{n-1}X)_n
\end{equation}
is a covering. Here $sk_n$ is the $n$-skeleton functor, i.e., $sk_n(X.)$ is the $n$-truncated
simplicial variety $(X_m)_{m\leq n}$. The functor
$$
cosk_n^X: \; (\;n\mbox{-truncated simplicial } X\mbox{-varieties}\;) \longrightarrow (\;\mbox{simplicial } X\mbox{-varieties}\;)
$$
is the right adjoint of $sk_n$ (which exists by general nonsense),
and the map \eqref{eq.cosk} comes from the adjunction map
$$
X. \longrightarrow cosk^X_{n-1}sk_{n-1} X.\;\;.
$$
It is a standard fact that the considered cohomology theories satisfy descent for hypercoverings.
So for a hypercovering $X. \rightarrow X$ as above the map
\begin{equation}\label{eq.hypercov}
H^n(X,\q) \mathop{\longrightarrow}\limits^{\sim} H^n(X.,\q)
\end{equation}
from the cohomology of $X$ to the cohomology of the simplicial complex variety $X.$ is
an isomorphism. Note that for any simplicial variety $X.$ we have a spectral sequence
\begin{equation}\label{eq.ss.hypercov}
E_1^{p,q} = H^q(X_p,\q) \Rightarrow H^{p+q}(X.,\q)\,.
\end{equation}

Call $X. \longrightarrow X$ a smooth hypercovering if all $X_p$ are smooth. A smooth hypercovering
exists, if every variety has a smooth covering, so this holds over any field by de Jong's
resolution of singularities.

\medskip
Now consider the case where $X$ is proper. Then all $X_n$ can be chosen to be smooth and projective,
in which case we call $X. \rightarrow X$ a smooth projective hypercovering. From \eqref{eq.hypercov}
and \eqref{eq.ss.hypercov} we obtain a spectral sequence
\begin{equation}\label{eq.weight.ss.2}
E_1^{p,q} = H^q(X_p,\q) \Rightarrow H^{p+q}(X,\q)\,.
\end{equation}
Note that $H^q(X_p,\q)$ carries a pure Hodge structure of weight $q$. Deligne defined the mixed Hodge
structure on $H^n(X,\q)$ in such a way that the spectral sequence \eqref{eq.weight.ss.2}
gives the weight filtration, and is a spectral sequence of mixed Hodge structures. Similarly as
in section 1 it follows that the spectral sequence degenerates at $E_2$.

\medskip
Analogous facts hold for the \'etale $\ql$-cohomology of varieties over finitely generated
fields and the weights on them.

\medskip
On the other hand it is known that the above descent theory does not extend to other functors -
like algebraic $K$-theory, or Chow groups, or motivic cohomology - unless one uses $\q$-coefficients.
But in \cite{GS} H. Gillet and C. Soul\'e developed a theory which works for integral coefficients.
For this they replaced coverings by so-called envelopes, i.e., surjective proper morphisms
$\pi: X' \rightarrow X$ of schemes such that every $x\in X$ has a point $x'\in X'$ mapping to $x$
such that the morphism $k(x) \rightarrow k(x')$ of residue fields is an isomorphism. In particular,
$\pi$ must be generically birational if $X$ and $X'$ are reduced. Call an envelope
$X' \rightarrow X$ of varieties smooth if $X'$ is smooth. If the ground field has characteristic
zero, then every reduced variety $X$ has a smooth envelope by Hironaka's resolution of singularities
\cite{Hi}. By a standard technique this also gives a smooth hyperenvelope of $X$, i.e.,
a simplicial $X$-scheme $X. \rightarrow X$
such that each
\begin{equation}\label{eq.cosk.2}
X_n \longrightarrow (cosk_{n-1}^Xsk_{n-1}X)_n
\end{equation}
is a smooth envelope. Now Gillet and Soul\'e showed that algebraic $K$-theory, Chow groups and many
related functors have descent for hyperenvelopes.

\medskip
Via these methods they were also able to construct the following. Let $Corr_K$ be the category of
correspondences over a field $K$: objects are smooth projective varieties over $K$ (not necessarily
geometrically connected), and morphisms are algebraic correspondences modulo rational equivalence.
Let $CHM_K$ be the idempotent completion of $Corr_K$ - objects are pairs $(X,p)$ where $X$ is a
smooth projective variety over $K$ and $p$ is an idempotent in $End_{Corr_K}(X)$. This is usually
called the category of Chow motives over $K$ (with integral coefficients); it is an additive category
where each idempotent has a kernel and cokernel. In contrast to \cite{GS} let us normalize the
categories so that the functor $X \mapsto M(X) = (X,id)$ from varieties to motives is covariant.
Let $K^b(CHM_K)$ be the homotopy category of bounded chain complexes in $CHM_K$. Then one has:

\begin{thm}(\cite{GS}) Let $K$ be a field of characteristic zero, or assume that resolution of
singularities exists over $K$.

\smallskip
(a) For any variety $X$ over $K$ there is an associated complex $W(X)$ in $K^b(CHM_K)$ (called the
{\it weight complex} of $X$), which is determined up to unique isomorphism.

\smallskip
(b) The association $X \mapsto W(X)$ is covariant for proper morphisms, i.e., a functor
on $V^\ast_K$.

\smallskip
(c) If $X$ is a variety, $Y \hookrightarrow X$ is a closed subvariety and $U \hookrightarrow X$ is
the open complement, then one has a canonical exact triangle
\begin{equation}\label{eq.ex.triangle}
W(Y) \longrightarrow W(X) \longrightarrow W(U) \longrightarrow W(Y)[1]\,.
\end{equation}

\smallskip
(d)If $U$ is a smooth variety of dimension $d$ and $U \subset X \supset Y$ is a good
compactification, then $W(U)$ is represented by the complex
$$
M(Y^{[d]}) \rightarrow M(Y^{[d-1]}) \rightarrow \ldots \rightarrow M(Y^{[1]}) \rightarrow M(X)\,,
$$
with the obvious differentials.
\end{thm}

Since the coefficients are integral everywhere, a remarkable consequence of this theory is that
one has a well-defined weight filtration on integral singular cohomology with compact supports $H^n_c(X,\z)$
(over $\mathbb C$) or \'etale cohomology $H^n_c(\overline X,\zl)$,
similarly for torsion coefficients, if one defines it via hyperenvelopes. The first one coincides
with Deligne's filtration after tensoring with $\q$, but it is shown in \cite{GS} that it cannot be
recovered from the $\q$-filtration. It would be interesting to see if this weight filtration - which
is trivial for smooth projective varieties by definition, gives some interesting information on
coefficients mod $\ell$, say.

\medskip
The results above also give the following, which immediately implies Theorem \ref{thm.homology.iso} (i).

\begin{thm}\label{thm.5.9} (\cite{GS} 3.1.1) If $H: CHM_K \longrightarrow Ab$ is a covariant functor from Chow motives
to abelian groups, there is a natural way to extend $H$ to a homology thery $H_.(-)$ on $V^\ast_K$
such that the following holds for {\it smooth projective} $X$:
$$
H_a(X) = \left\{\begin{array}{ccc} 0 & , & a\neq 0,\\H(X)& , & a=0\end{array}\right.\, .
$$
\end{thm}

\smallskip
In fact, one gets the weight homology $H^W_.(-)$ by applying \ref{thm.5.9} to the functor $H(X) = (\qzl)^{\pi_0(X)}$.

\medskip
As for the morphism of homology theories in Theorem \ref{thm.homology.iso} (ii), it is obtained
by refining Theorem \ref{thm.5.9} to a functor with values in complexes, applying it to the complexes
$\Cb(X,\qzl)$, and defining $\varphi$ as induced by functoriality of the construction starting
from the trace map
$$
f_\ast: \Cb(X,\qzl) \rightarrow \Cb(\Spec(K),\qzl) = \qzl
$$
for a connected smooth projective variety $f: X \rightarrow \Spec(K)$.

\section{Varieties over finite fields}

Kato also stated a conjecture for varieties over finite fields. For such varieties $X$ he
defined a complex $C^{1,0}(X,\qzl)$:
\begin{equation}\label{eq.Kato.cx.2}
\begin{array}{ccccl}
\ldots&  \rightarrow& \mathop\oplus\limits_{x\in X_a} H^{a+1}(k(x),\mathbb Q_\ell/\mathbb Z_\ell(a))& \rightarrow&\mathop\oplus\limits_{x\in X_{a-1}} H^{a}(k(x), \mathbb Q_\ell/\mathbb Z_\ell(a-1))\\
\rightarrow& \ldots & \ldots & \rightarrow &\mathop\oplus\limits_{x\in X_0} H^1(k(x),\qzl)\, .
\end{array}
\end{equation}
and stated the following conjecture.

\begin{conj}\label{conj.Kato.finite.field}(\cite{Ka1} Conj. 0.3)
If $X$ is connected, smooth and proper over a finite field $k$, then
\begin{equation}\label{eq.Hasse.X.2}
H_a(C^{1,0}(X,\mathbb Q_\ell/\mathbb Z_\ell))=\left\{\begin{array}{ccc} 0 & , & a>0,\\\mathbb Q_\ell/\mathbb Z_\ell& , & a=0\end{array}\right.\, .
\end{equation}
\end{conj}

For $\dim(X) = 1$ this conjecture amounts to \eqref{eq.Brauer} with $K = k(X)$, for $\dim(X)=2$
the conjecture was proved by Colliot-Th\'el\`ene, Sansuc, and Soul\'e \cite{CTSS} for $\ell$ invertible in $k$,
and by M. Gros \cite{Gr} and K. Kato \cite{Ka1} if $\ell = \ch(k)$. S. Saito \cite{Sa} proved that
$H_a(C^{2,1}(X,\qzl)) = 0$ for $\dim(X)=3$ and $\ell \neq \ch(k)$. For $X$ of any dimension
Colliot-Th\'el\`ene \cite{CT} (for $\ell\neq\ch(k)$) and Suwa \cite{Su}
(for $\ell=\ch(k))$ proved that $H_a(C^{1,0}(X,\qzl))=0$ for $0<a\leq 3$.

\medskip
In \cite{Ja4} Thm. 6.1 it is shown that resolution of singularities for varieties of dimension $\leq d$
would imply this conjecture for $X$ smooth projective of dimension $\leq d$.
In \cite{JS2} some recent results on resolution of singularities \cite{CJS} are applied in a
different way to obtain the following unconditional result.

\begin{thm}\label{thm.Hasse.JS}
If $X$ is connected, smooth and projective over a finite field, then
$$
H_a(C^{1,0}(X,\mathbb Q_\ell/\mathbb Z_\ell))=\left\{\begin{array}{ccc} 0 & , & 0<a\leq 4,\\
\mathbb Q_\ell/\mathbb Z_\ell& , & a=0\end{array}\right.\, .
$$
\end{thm}

There are applications to the finiteness of certain motivic cohomology groups with finite
coefficients of $X$:

\begin{thm}\label{thm.motivic.coh}(\cite{JS2} Thm. 6.3) Let $X$ be a smooth projective
variety of pure dimension $d$ over a finite field $k$. Assume that the Galois symbol
\begin{equation}\label{eq.Galois.symbol}
K^M_{q}(L) \longrightarrow H^{q}(L,\z/\ell\z(q))
\end{equation}
between Milnor $K$-theory and Galois cohomology (\cite{Mi}, \cite{Ta}, \cite{BK}) is
surjective for all $\ell \mid n$ and all fields $L$ above $k$. Then the cycle maps
$$
\rho_X^{r,t}: CH^r(X,t;\mathbb Z/n\mathbb Z) = H^{2r-t}_{\mathcal M}(X,\mathbb Z/n\mathbb Z(r)) \longrightarrow H^{2r-t}(X,\mathbb Z/n\mathbb Z(r))
$$
between higher Chow groups/motivic cohomology groups with finite coefficients and \'etale cohomology are isomorphisms
for $r>d$ and $t\leq q-1$, and for $r=d$ and $t \leq q-2\leq 2$. In particular the above higher Chow groups
are finite under these conditions.
\end{thm}

See also \cite{Ja3} for some results in similar direction for $\z$-coefficients, and work of Th. Geisser \cite{Ge},
who formulated and studied an integral form of Kato's conjecture (over a finite field).

\medskip
The surjectivity of \eqref{eq.Galois.symbol} is known for $\ell = \ch(k)$ (\cite{BK}), and for
$\ell$ invertible in $k$ in the following cases: $q=1$ (Kummer theory), $q=2$ (\cite{MS}), $\ell = 2$ (\cite{V1});
it has been announced by Rost and Voevodsky to hold in general \cite{Ro},\cite{V2}, see also \cite{SJ} and \cite{Weib}.

\medskip
Finally let me mention that Kato \cite{Ka1} also stated some conjectures for regular proper schemes
over $\mathbb Z$ or $\mathbb Z_p$ (related to those considered in sections 7, 8 and 9).
These were studied in \cite{JS}, again by weight methods.

\medskip
\begin{footnotesize}I thank Takeshi Saito heartily for carefully checking the manuscript.
\end{footnotesize}

\vspace{2cm}

Uwe Jannsen\\
NWF I - Mathematik\\
Universit\"at Regensburg\\
93040 Regensburg\\
GERMANY\\
uwe.jannsen@mathematik.uni-regensburg.de

\end{document}